\documentclass{www-notes}
\usepackage{fullpage}
\usepackage{upgreek}

\newcommand*{\mr}{\mathring}
\newcommand*{\paru}{\eth}
\newcommand*{\parub}{\ub{\eth}}
\newcommand*{\Osub}[3]{\mathcal{O}_{#1,#2}(#3)}

\makeatletter
\newcommand{\subjclass}[2][2010]{%
  \let\@oldtitle\@title%
  \gdef\@title{\@oldtitle\footnotetext{#1 \emph{Mathematics subject classification.} #2}}%
}

\title{Geometric analysis of 1+1 dimensional quasilinear wave equations}
\author{Leonardo Enrique Abbrescia\Affiliation{Department of Mathematics, Vanderbilt University, Nashville, Tennessee, USA; \url{leonardo.abbrescia@vanderbilt.edu}} \and Willie Wai Yeung Wong\Affiliation{Department of Mathematics, Michigan State University, East Lansing, Michigan, USA; \url{wongwwy@math.msu.edu}}}

\begin{document}


\maketitle

\begin{wwwabstract}
	We prove global well-posedness of the initial value problem for a class of variational quasilinear wave equations, in one spatial dimension, with initial data that is not necessarily small. Key to our argument is a form of quasilinear null condition (a ``nilpotent structure'') that persists for our class of equations even in the large data setting. This in particular allows us to prove global well-posedness for $C^2$ initial data of moderate decrease, provided the data is sufficiently close to that which generates a simple traveling wave.
	We take here a geometric approach inspired by works in mathematical relativity and recent works on shock formation for fluid systems. First we recast the equations of motion in terms of a dynamical double-null coordinate system; we show that this formulation semilinearizes our system and decouples the wave variables from the null structure equations. After solving for the wave variables in the double-null coordinate system, we next analyze the null structure equations, using the wave variables as input, to show that the dynamical coordinates are $C^1$ regular and covers the entire spacetime. 
\end{wwwabstract}

\section{Introduction}

This paper concerns quasilinear wave equations on $\Real^{1,1}$. Denoting by $\eta$ the Minkowski metric and by $\sigma = \sigma(\phi) = (\eta^{-1})(\mathrm d\phi,\mathrm d\phi)$, we are specifically interested in solutions to the equation
\begin{equation}\label{eq:mainwave}
	\partial_{\mu}\left( e^{f(\sigma)} \eta^{\mu\nu} \partial_\nu \phi \right) = 0.
\end{equation}
This equation arises as the Euler-Lagrange equation for the action
\[ \int_{\Real^{1,1}} F(\sigma) ~\mathrm{dvol}_\eta \]
where $F\colon\Real\to\Real$ is a primitive of $e^f$.  

Equation \eqref{eq:mainwave} is a particular example of a generic quasilinear equation on $\Real^{1,1}$ of the form 
\begin{equation}\label{eq:genqnlw}
	\partial_{\mu} \left( \tilde{f}(\partial\phi) \eta^{\mu\nu} \partial_\nu \phi \right) = 0.
\end{equation}
Equations of the form \eqref{eq:genqnlw} fall under the umbrella of \emph{systems of one-dimensional conservation laws}, which is an extensively studied field whose roots date back to Riemann \cite{bR1860}. The authors refer the reader to the monographs \cite{dafermos2005hyperbolic, li1994global} and the references therein for a comprehensive (but not exhaustive) overview of the literature. 
 
 It is well-known that generic quasilinear wave equations on $\Real^{1,1}$ of the form \eqref{eq:genqnlw} do not admit global-in-time solutions arising from non-trivial compactly-supported initial data
(see \cite{John1974, MajTho1987, SpHoLW2016} and references therein). Genericity here refers to the ``genuinely nonlinear'' condition of Lax \cite{Lax1964, Lax1973}; solutions to such equations with compactly-supported initial data always terminate in finite-time shock singularities.
Our equation, on the other hand, is not genuinely nonlinear; rather by virtue of the $\sigma$ dependence in the nonlinearity, it satisfies the ``null condition'' of Klainerman. 
The main result of the present work is that \eqref{eq:mainwave} admits an open set
of global-in-time solutions, corresponding to perturbations of a simple wave ansatz. Global Lipshitz wellposedness in the small initial data setting is a corollary of our results when the simple wave identically vanishes. In fact, we are able to show the following (see \thref{thm:main} for a precise statement):

\begin{thm}[Rough statement]
	The initial value problem for \eqref{eq:mainwave} is globally well-posed, provided the initial data is of \emph{moderate decrease} and is sufficiently close (in $C^2$ norm) to that of a simple traveling wave. 
\end{thm}

We note that the idea that small data global existence results can be extended to (on occasions, \emph{semi-}global) existence results for data sufficiently close to ``simple wave'' initial data is not entirely new, especially for the case of higher dimensional wave equations. 
We mention specifically the results of \cite{abbrescia2019global,Sideri1989} in the context of the semilinear wave maps equation, where the notion of ``simple wave'' is one with one-dimensional image, and the results \cite{MiPeYu2019, WanYu2016} of more general semilinear equations where the notion of ``simple waves'' is of an out-going, dispersive, nearly spherically-symmetric solution. In the present work by ``simple traveling waves'' or just ``simple waves'' we mean a function on the $(1+1)$-dimensional Minkowski space
\[\mathring\phi(t,x) \eqdef \zeta(t-x),\]
where $\zeta \colon \Real \to \Real$ is at least $C^2$.

The role played by the null condition in our setting is intricate: unlike the cases of small data (and the aforementioned perturbations of simple wave) theory for higher dimensional nonlinear wave equations with null condition, there are no uniform-in-time decay estimates for solutions to the linear wave equation in one spatial dimension. 
The global existence mechanism therefore relies on the essential fact that two wave packets with distinct velocities can only interact for a finite time period.
This idea has been captured previously in special cases: \textbf{(I)} certain one-dimensional conservation laws which are not genuinely nonlinear \cite{LiuJianli2007Abog}; \textbf{(II)} the situation for the small data global existence for the so-called membrane equation on $\Real^{1,1}$ was first analyzed in  \cite{Lindbl2004}, whose results were generalized by the second author in \cite{Wong2016p}; \textbf{(III)} recently a study of the small data global existence for \emph{semilinear} wave equations obeying the null condition using conformally weighted energy estimates was undertaken in \cite{LuYaYu2017}.
A key argument in the latter work is to demonstrate the strong localization, around sets that are largely disjoint in spacetime, of the forward- and backward-travelling components of the solution. 
The arguments of \cite{LuYaYu2017} have also been recently extended to the small data quasilinear setting by Zha \cite{Zha2019}. The aim of the present manuscript is to systematically analyze perturbations of simple (not necessarily small) traveling wave solutions of general quasilinear equations of the form \eqref{eq:mainwave} using this idea.

In order to capture this ``non-interaction'', we use an approach based on the geometric analysis of the principal symbol corresponding to \eqref{eq:mainwave}. 
This approach is the same underlying the recent works on higher dimensional shock formation for quasilinear waves 
\cite{Christ2007a, Speck2016, HoKlSW2016, SpHoLW2016}; our results complement them in that we study the case where genuine nonlinearity fails, which is explicitly excluded in those works by assumption. We also replace the $L^2$ type estimates of \cite{LuYaYu2017} with more direct $L^\infty$ type estimates available to us in the $(1+1)$-dimensional setting. This geometric approach was previously taken by the second author to treat the membrane equation \cite{Wong2016p}, which corresponds to \eqref{eq:mainwave} with the primitive $F(\sigma) = \sqrt{1 + \sigma}$. 
There, however, the solution has a clear interpretation as an immersed submanifold of a flat, higher-dimensional Minkowski space, and the relevant equations of motion can be obtained rather straightforwardly from purely geometric considerations. 
Furthermore, the geometry also forces stronger, algebraic decoupling of the various components of the system, a fact which can also be regarded as a specific manifestation of the extremely strong null condition enjoyed by the membrane equation.  In the more general setting treated by the present manuscript, one of the key steps is a geometric formulation of quasilinear equations of the form \eqref{eq:mainwave} that is compatible with $L^\infty$ types estimates and exposes the weaker (but nevertheless present) non-resonance of the null condition. 

A further complication is introduced by the fact that we work in the setting not of small initial data, but we allow initial data to be a small perturbation of a (potentially large) simple wave. The perturbed system (see \eqref{eq:semisys:p}) now has linear terms with large coefficients. In view of the non-decay of solutions to the wave equation on $\Real^{1,1}$, this can potentially cause instability, or at least make it not feasible to address the equation perturbatively. In the case of the membrane equation we recently showed in \cite{AbbWong2019} that these linear terms vanish in the perturbation equation when written in an appropriate gauge. There, the large traveling wave appears only as a coefficient of  \emph{quadratic and higher order} nonlinearities and, although there is some growth coming from these coefficients, the strong null condition of the membrane equation provides just enough decay in three spatial dimensions for $C^2$ stability of simple traveling waves. For the system studied in the present work, such linear terms are unavoidable. However, we exhibit a \emph{nilpotent} (or \emph{weak null}) structure. This allows us to design an appropriate iteration scheme to treat even these large coefficient terms as ``small" perturbations to the dynamics. 

\textbf{Acknowledgements--} LA gratefully acknowledges support from an NSF Postdoctoral Fellowship. WWY Wong is supported by a  Collaboration Grant from the Simons Foundation, \#585199.

\section{The acoustic metric and its geometry}\label{sec:geometry}
It is convenient to rewrite \eqref{eq:mainwave}; dividing by $e^{f(\sigma)}$ we can rewrite the equation as
\begin{equation}\label{eq:mainwaveH}
	g^{\mu\nu} \partial^2_{\mu\nu} \phi = 0,
\end{equation}
where the \emph{inverse acoustic metric} is given by
\begin{equation}\label{eq:dfn:acinvmetric}
	g^{\mu\nu} \eqdef \eta^{\mu\nu} + 2 f'(\sigma)\; \eta^{\mu\alpha}\partial_\alpha \phi \;\eta^{\nu\beta} \partial_\beta \phi.
\end{equation} 
A standard computation yields that the corresponding \emph{acoustic metric}, which satisfies $g_{\lambda\mu} g^{\mu\nu} = \delta_\lambda^\nu$, is given by
\begin{equation}\label{eq:dfn:acmetric}
	g_{\mu\nu} \eqdef  \eta_{\mu\nu} -\frac{2 f'(\sigma)}{1 + 2 f'(\sigma) \sigma} \partial_\mu\phi \partial_\nu\phi.
\end{equation}
It is worth remarking that 
\begin{gather}
	g^{\mu\nu} \partial_\mu\phi \partial_\nu \phi = \sigma + 2 f'(\sigma) \sigma^2, \label{eq:gcontract} \\
	\det g^{-1} = -1 - 2 f'(\sigma) \sigma, \\
	\det g = - (1 + 2 f'(\sigma) \sigma)^{-1}. 
\end{gather}
It is convenient to introduce the auxiliary vector-valued variable $\Phi$ with
\begin{equation}
	\Phi_0 = \partial_t \phi, \qquad \Phi_1 = \partial_x \phi.
\end{equation}
Note that $\sigma = - \Phi_0^2 + \Phi_1^2$. Taking a derivative of \eqref{eq:mainwave} we see that $\Phi$ satisfies
\[
	\partial_\mu \left( e^f g^{\mu\nu} \partial_\nu \Phi \right) = 0.
\]
\begin{exa}
	To illustrate, in the case of the membrane equation, we have 
	\[ f(\sigma) = - \frac12 \ln(1 + \sigma) \]
	so the quantity $1 + 2 f'(\sigma) \sigma = (1+\sigma)^{-1}$,
	and thus the acoustic metric reduces to $\eta_{\mu\nu} + \partial_\mu \phi \partial_\nu\phi$.
	One can further check that in this case 
	\[  e^f \eta^{\mu\nu} \partial_\nu \phi = \sqrt{\abs{\det g}} g^{\mu\nu} \partial_\nu \phi \]
	and hence the \eqref{eq:mainwave} can also be re-written as 
	\[ \Box_{g(\partial\phi)} \phi = 0 \]
	where $\Box_{g}$ is the geometric wave operator for the acoustic metric $g$. 
	This expression is fundamentally what gives the exceptional structure used in \cite{Wong2016p}. 
\end{exa}
Reorganizing this equation, we have
\begin{align*}
	0 & = e^{-f} \partial_{\mu} \left( e^f g^{\mu\nu} \partial_\nu \Phi \right) \\
	 &= \partial_\mu \left( g^{\mu\nu} \partial_\nu \Phi \right) + f' g^{\mu\nu} \partial_\mu \sigma \partial_\nu \Phi \\
	 & = \sqrt{\abs{\det g^{-1}}} \partial_\mu \left( \sqrt{\abs{\det g}} g^{\mu\nu} \partial_\nu \Phi \right) 
	 	+  \frac12 g^{\mu\nu} \partial_\mu (\ln \abs{\det g^{-1}}) \partial_\nu \Phi
		+ f' g^{\mu\nu} \partial_\mu\sigma \partial_\nu\Phi 
\end{align*}
or, that the geometric wave equation 
\begin{equation}\label{eq:derwaveR}
	 \Box_g \Phi + \underbrace{\left[ \frac{f''(\sigma) \sigma + f'(\sigma)}{1 + 2 f'(\sigma) \sigma} + f'(\sigma)\right]}_{G(\sigma)} g^{\mu\nu} \partial_\mu\sigma \partial_\nu \Phi = 0
\end{equation}
is satisfied. We note that for fixed $\sigma$, this equation is linear in $\Phi$, and hence taking linear combinations we can define the scalars
\begin{equation}
	\Psi \eqdef \Phi_0 + \Phi_1, \qquad \ub{\Psi} \eqdef \Phi_0 - \Phi_1
\end{equation}
and have that they satisfy
\begin{equation}\label{eq:derwaveDN}
\begin{gathered}
	\Box_g \Psi + G(\sigma) g^{-1}(\D*\Psi, \D*\sigma) = 0; \\
	\Box_g \ub{\Psi} + G(\sigma) g^{-1}(\D*\ub\Psi, \D*\sigma) = 0.
\end{gathered}
\end{equation}
On the other hand, since we can re-write 
\begin{equation}\label{eq:sigmacomp}
	\sigma = - \Psi \ub{\Psi}, 
\end{equation}
it satisfies
\begin{equation}\label{eq:swaveDN}
\begin{aligned}
 \Box_g \sigma &= - \Psi \Box_g \ub{\Psi} - \ub{\Psi} \Box_g \psi - 2 g^{-1}(\D*\Psi, \D*\ub{\Psi})\\
 	&= G(\sigma) g^{-1}(\Psi \D\ub{\Psi} + \ub{\Psi} \D\Psi, \D*\sigma) - 2 g^{-1}(\D*\Psi, \D*\ub{\Psi})\\
	&= - G(\sigma) g^{-1}(\D*\sigma, \D*\sigma) - 2 g^{-1}(\D*\Psi, \D*\ub{\Psi}).
\end{aligned}
\end{equation}

Our approach to understanding the long-time behavior for \eqref{eq:mainwave} goes through the study of \eqref{eq:derwaveDN} and \eqref{eq:sigmacomp}; we could alternatively also use \eqref{eq:swaveDN}. 
At this juncture we will take advantage of the \emph{conformal invariance} of the Laplace-Beltrami operator in two spacetime dimensions. 
This will allow us to effectively \emph{semilinearize} equations \eqref{eq:derwaveDN} and \eqref{eq:swaveDN} and separate their analysis from the analysis of the spacetime geometry. 

More precisely, we will rewrite the equations of motion \eqref{eq:derwaveDN} in \emph{dynamical double null coordinates} generated by the conformal structure of the acoustic metric $g$. 
We will see that in this coordinate system the metric decouples from the equations for $\Psi, \ub{\Psi}$, and $\sigma$, which themselves form an autonomous system of semilinear wave equations which can be solved independently of the metric $g$. 
The solution is then completed by studying the transition map relating the dynamical double null coordinates to the standard rectangular coordinates of $\Real^{1,1}$. 
We remark that this approach is similar to the approach pioneered in the proof of the nonlinear stability of Minkowski space in general relativity \cite{ChrKla1993} as well as in recent works establishing stable shock formation for quasilinear wave equations (e.g.\ \cite{Christ2007a, HoKlSW2016, SpHoLW2016}). 

Denote now by $u, \ub{u}$ two independent scalar functions satisfying $g^{-1}(\D*u, \D*u) = g^{-1}(\D*\ub{u}, \D*\ub{u}) = 0$. 
The functions are defined up to reparametrization $u \mapsto v(u), \ub{u} \mapsto \ub{v}(\ub{u})$, which leaves their level sets invariant; we will later make use of this freedom to normalize them by setting their values at $t = 0$. 
Together $u$ and $\ub{u}$ define a \emph{double null} coordinate system; they are \emph{dynamic} in the sense that their definition depends on the acoustic metric $g$ which, itself, depends on the unknowns $\Psi, \ub{\Psi}$, and  $\sigma$. 
Relative to the double null coordinates, the metric $g$ takes the form 
\begin{equation}
	\begin{gathered}
		g = \Omega (\D*u \otimes \D*\ub{u} + \D*\ub{u} \otimes \D*{u}); \\
		g^{-1} = \Omega^{-1} (\paru \otimes \parub + \parub \otimes \paru).
	\end{gathered}
\end{equation}
Here and throughout, we will use $\paru$ to refer to the vector field corresponding to coordinate partial differentiation with respect to $u$ (holding $\ub{u}$ fixed), and $\parub$ \emph{vice versa}; in particular the two fields commute $[\paru, \parub] = 0$. 
The conformal invariance of the Laplace-Beltrami operator implies that
\[ \Box_g = 2\Omega^{-1} \paru \parub \]
and hence we see that the equations of motion reduce to the following semilinear system in $u, \ub{u}$ coordinates:
\begin{equation}\label{eq:semisys}
	\begin{gathered}
		2\paru\parub \Psi + G(\sigma) \; \paru \sigma \; \parub \Psi + G(\sigma) \; \paru \Psi \; \parub\sigma = 0,\\
		2\paru\parub \ub\Psi + G(\sigma) \; \paru \sigma \; \parub \ub\Psi + G(\sigma) \; \paru \ub\Psi \; \parub\sigma = 0, \\
		\sigma = - \Psi \ub{\Psi}.
	\end{gathered}
\end{equation}
After solving for $\Psi, \ub{\Psi}$, and $\sigma$ in terms of $u$ and $\ub{u}$ through \eqref{eq:semisys}, we can recover the original unknown $\phi$ as follows. Observe that the Hessian with regard to the acoustic metric can be written as
\[ \covD^2_{\mu\nu} \phi = \partial^2_{\mu\nu} \phi - \Gamma^{\lambda}_{\mu\nu} \partial_\lambda \phi,\]
where the Christoffel symbol is defined by 
\[ \Gamma^{\alpha}_{\beta\gamma} = \frac12 g^{\alpha \delta} \left[ \partial_\beta g_{\delta\gamma} + \partial_\gamma g_{\delta\beta} - \partial_{\delta} g_{\beta\gamma} \right]. \]
Writing 
\begin{equation}\label{eq:dfn:H}
	H(\sigma) \eqdef - \frac{2 f'(\sigma)}{1 + 2 f'(\sigma) \sigma} 
\end{equation}
we have, by virtue of \eqref{eq:dfn:acmetric} that
\begin{multline}\label{eq:dfn:chrsym}
	\Gamma^{\alpha}_{\beta\gamma} = \frac12 g^{\alpha \delta} \Big[ H'(\sigma) \left( \partial_\beta \sigma \Phi_\delta \Phi_\gamma + \partial_\gamma \sigma \Phi_\delta \Phi_\beta - \partial_\delta \sigma \Phi_\beta \Phi_\gamma \right) \\
		+ H(\sigma) \Phi_\delta \left( \partial_\beta \Phi_\gamma + \partial_\gamma \Phi_\beta \right) \\
		+ H(\sigma) \Phi_\gamma \left( \partial_\beta \Phi_\delta - \partial_\delta \Phi_\beta \right) 
		+ H(\sigma) \Phi_\beta \left (\partial_\gamma \Phi_\delta - \partial_\delta \Phi_\gamma \right) \Big].
\end{multline}
From the integrability condition $\partial_\mu\Phi_\nu = \partial_\nu \Phi_\mu$ (partial derivatives commute) we see that the terms on the final line above vanish. 
Therefore we have that 
\[ \Box_g \phi = g^{\mu\nu} \partial^2_{\mu\nu} \phi - g^{\mu\nu} \Gamma^\lambda_{\mu\nu} \partial_{\lambda} \phi.\]
Through \eqref{eq:mainwaveH} we see $g^{\mu\nu}\partial^2_{\mu\nu} \phi = g^{\mu\nu} \partial_\mu \Phi_\nu = 0$, and hence $\phi$ solves the wave equation
\[ \Box_g \phi + \frac12 g^{\mu\nu} g^{\lambda\delta} H'(\sigma) \Phi_{\lambda} \Phi_{\delta} \partial_\mu \phi \partial_\nu \sigma = 0, \]
which we can further simplify using \eqref{eq:gcontract} to 
\begin{equation}
	\Box_g \phi + \frac12 \sigma (1 + 2 f'(\sigma)\sigma) H'(\sigma) g^{-1}(\D*\sigma, \D*\phi) = 0
\end{equation}
or, in terms of the null coordinates, 
\begin{equation}\label{eq:semisysphi}
	\paru\parub \phi + \frac14 \sigma(1 + 2 f'(\sigma) \sigma) H'(\sigma) \left[ \paru \sigma \parub \phi + \paru \phi \parub \sigma \right] = 0.
\end{equation}
We note that with $\Psi$, $\ub{\Psi}$, and $\sigma$ being known quantities from \eqref{eq:semisys}, the equation \eqref{eq:semisysphi} is \emph{linear} in $\phi$. 

Solving \eqref{eq:semisys} and \eqref{eq:semisysphi} would give us $\phi, \Psi, \ub{\Psi}$, and $\sigma$ as \emph{functions of the coordinates $u, \ub{u}$}.
To recover them as functions of the original $t,x$ coordinates of $\Real^{1,1}$, we need to examine the corresponding change of variables map. 

\begin{rmk}
	We note that from their definitions we have that $\Phi_0 = \frac12 (\Psi + \ub\Psi)$ and $\Phi_1 = \frac12 (\Psi - \ub\Psi)$ are also recovered as functions of $u, \ub{u}$ once the system \eqref{eq:semisys} is solved. By virtue of their definitions \eqref{eq:dfn:acinvmetric} and \eqref{eq:dfn:acmetric}, this implies that the \emph{rectangular components} of the acoustic metric and its inverse can also be recovered as functions of $u, \ub{u}$. 
	
	In particular, the acoustic metric is non-degenerate provided $1 + 2 f'(\sigma) \sigma > 0$, and that $\Psi, \ub{\Psi}, \sigma$ are bounded. 
\end{rmk}

To control the change of variables, it suffices to control the Jacobian quantities $\paru t, \paru x, \parub t, \parub x$. 
With these quantities controlled one can in principle invert the coordinate transformation and recover $u, \ub{u}$ as functions of $t, x$. 
The obstacle to this is primarily \emph{singularities where the Jacobian determinant vanishes}; this corresponds to \emph{shock formation} for the system. 

Define now the (acoustic) null vector fields $L$ and $\ub{L}$ whose components are given by
\begin{equation}\label{eq:dfn:nullvec}
	L^{\mu} = g^{\mu\nu} \partial_\nu u, \qquad \ub{L}^\mu = g^{\mu\nu} \partial_\nu \ub{u}. 
\end{equation}
As $L$ is the gradient of a solution to an eikonal equation, it is geodesic with respect to $g$, similarly $\ub{L}$. 
Notice that 
\[ L^\lambda = L(x^\lambda) = g^{\mu\nu} \partial_\nu u \partial_\mu x^\lambda = g^{-1}(\D*u, \paru x^\lambda \D u + \parub x^\lambda \D\ub{u}) = \Omega^{-1} \parub x^\lambda \]
and further
\[ g_{\mu\nu} L^\mu \ub{L}^\nu = g( L, \ub{L}) = g^{-1}(\D*u, \D*\ub{u}) = \Omega^{-1}.\]
Hence control of the coordinate components $L^\mu$ and $\ub{L}^\mu$ will allow us to further control not only the Jacobian values $\paru t, \parub t, \paru x, \parub x$, but also the conformal factor $\Omega$. 

To control the components $L^\mu$, we will use the approach of \cite{Christ2007a}, taking advantage of the fact that $L$ is geodesic. Denote by $\Gamma^\alpha_{\beta\gamma}$ the Christoffel symbols of the acoustic metric $g$ in the rectangular coordinate system, we have $L^\mu$ must satisfy
\[ \parub L^\mu + \Omega \Gamma^{\mu}_{\nu\lambda} L^\nu L^\lambda = 0.\]
We rewrite the geodesic equation in the following form, using the expression \eqref{eq:dfn:chrsym} for the Christoffel symbol:
\[
	\parub L^\mu + g^{\mu\delta} \Phi_\delta L^\gamma \parub[ H(\sigma) \Phi_\gamma] 
		=  \frac12 \Omega g^{\mu\delta} \partial_\delta \sigma H'(\sigma) (L^\beta \Phi_\beta)^2.
\]
Simplifying further, and swapping $L$ for $\ub{L}$, we arrive at the system
\begin{equation}\label{eq:geodesic}
\begin{gathered}
	\parub L^\mu + g^{\mu\delta} \Phi_\delta L^\gamma \parub[ H(\sigma) \Phi_\gamma]
		=  \frac{H'(\sigma) (L^\beta \Phi_\beta)^2}{2g_{\alpha\gamma} L^\alpha \ub{L}^\gamma} (L^\mu \paru \sigma + \ub{L}^\mu \parub \sigma);\\
	\paru \ub{L}^\mu + g^{\mu\delta} \Phi_\delta \ub{L}^\gamma \paru[ H(\sigma) \Phi_\gamma]
		=  \frac{H'(\sigma) (\ub{L}^\beta \Phi_\beta)^2}{2g_{\alpha\gamma} L^\alpha \ub{L}^\gamma} (L^\mu \paru \sigma + \ub{L}^\mu \parub \sigma).
\end{gathered}
\end{equation}
We remark that all the terms appearing in the expression, except for the values of $L$ and $\ub{L}$, are known quantities that can be computed from $\sigma, \Psi, \ub{\Psi}$, once \eqref{eq:semisys} is solved.

\begin{rmk}\label{rmk:cov:reg}
	To guarantee that the change of variables map is regular, it suffices that $u, \ub{u}$ remain $C^1$ with non-vanishing Jacobian determinant. This requires $\Omega L^\mu$ and $\Omega \ub{L}^\mu$ to remain bounded, with $\Omega L$ and $\Omega \ub{L}$ to be linearly independent. 

	A sufficient condition for this is for the components $L^\mu, \ub{L}^\mu$ to remain bounded; $\Omega$ to remain bounded above and below, and $L$ and $\ub{L}$ to be linearly independent. 
	Now, since $L$ and $\ub{L}$ are null vectors, provided the acoustic metric coefficients $g_{\mu\nu}$ remain bounded, their linear independence would be implied by the non-vanishing of their $g$-inner product. 

	To summarize, a sufficient condition for the regularity of the change of variables map is that the components $L^\mu, \ub{L}^\mu$, and $g_{\mu\nu}$ remain bounded, and that $\Omega$ is bounded above and below. 

	Additionally, to guarantee the map is $(u, \ub{u}) \to (t, x)$ is surjective, it further suffices that $L^0$ and $\ub{L}^0$ are both bounded away from $0$. 
\end{rmk}

\begin{assm}\label{assm:f:assm}
	We assume that $f$ is at least thrice continuously differentiable. This implies that for all $M_0$ sufficiently large (compared to $f'(0)$, $f''(0)$, and $f'''(0)$) there exists $m_0 > 0$ sufficiently small such that whenever $\abs{\sigma} \leq m_0$, the uniform bounds
\[\abs{G(\sigma)}, \abs{H(\sigma)}, \abs{f'(\sigma)}, \abs{G'(\sigma)}, \abs{H'(\sigma)}, \abs{1 + 2 f'(\sigma)\sigma}, \abs{1 + 2 f'(\sigma)\sigma}^{-1}, \abs{f''(\sigma)}, \abs{H''(\sigma)} \leq M_0 \]
	hold. 
\end{assm}

\section{Simple wave solutions}

Let us consider a special class of solutions corresponding to simple waves, that is to say, solutions of the form $\mathring\phi(t,x) = \zeta(t-x)$ (equivalently $t+x$) for some real valued function $\zeta$ which describes the profile of the wave. In this section and in the rest of the manuscript, quantities introduced in Section \ref{sec:geometry} associated to the simple waves will be adorned with a ring. For example, 
\[\mathring \Phi_0 = \partial_t\mathring\phi, \quad \mathring\sigma = \eta^{-1}(\mathrm d \mathring\phi,\mathrm d\mathring \phi),\]
and so forth.  By inspection it is clear that any traveling wave $\mathring\phi(t,x)$ solves \eqref{eq:mainwave}.
In this section we describe the simple waves in terms of the variables introduced in the previous section. 

First we see that $\mr\Phi_0(t,x) = \zeta'(t-x)$ and $\mr\Phi_1(t,x) = -\zeta'(t-x)$, and hence
\begin{equation}
	\mr \Psi \equiv 0, \quad \ub{\mr{\Psi}}(t,x) = 2\zeta'(t-x), \quad \text{ and } \mr\sigma \equiv 0.
\end{equation}
Denote by $\ub{u} = t - x$ and $v = t + x$, we note that the acoustic metric can be written as
\begin{equation}
	- \frac12 (\D*\ub{u} \otimes \D*v + \D*v \otimes \D*\ub{u}) + H(0) (\zeta'(\ub{u}))^2 \D*\ub{u} \otimes \D*\ub{u}
\end{equation}
where $H$ is defined by \eqref{eq:dfn:H}. Now if we let
\begin{equation}
	Z(\ub{u}) = - \int_0^{\ub{u}}  H(0) (\zeta'(s))^2 \D{s} 
\end{equation}
then we can factor the acoustic metric as
\begin{equation}
	\mr g = - \frac12 (\D*u \otimes \D*\ub{u} + \D*\ub{u} \otimes \D*u)
\end{equation}
where 
\begin{equation}
	u = v + Z(\ub{u}).
\end{equation}
A consequence is that $\mr \Omega \equiv - \frac12$ in this setting.

Next we can compute the null vector fields using \eqref{eq:dfn:nullvec}. Noting that $2f'(0) = - H(0)$, we obtain
\begin{equation}\label{eq:Lbg}\begin{aligned}
	\ub{\mr{L}}^0 &= -1, &\ub{\mr{L}}^1 &= -1;\\
	\mr L^0 &= -1 - H(0) (\zeta'(\ub{u}))^2, & \mr L^1 &= 1 - H(0) (\zeta'(\ub{u}))^2.
\end{aligned}\end{equation}

\begin{rmk}\label{rmk:hyp}
	From the computations above, we see that the level sets of the coordinate function $t$ are space-like with respect to the acoustic metric $g$, if and only if $H(0)(\zeta'(\ub{u}))^2 > -1$. In particular, when the function $f$ of \eqref{eq:mainwave} satisfies $f'(0)\leq 0$ this holds regardless of the profile $\zeta$; while for $f'(0) > 0$ the $\zeta$ is required to have sufficiently small slope to guarantee that the level sets of $t$ is space-like, and therefore the hyperbolicity of the perturbation equations relative to the $t$ level sets.  

	On the other hand, we also see that for any plane-wave solution, the Jacobian matrix
	\[ \begin{pmatrix}
		\paru t & \parub t\\ \paru x & \parub x 
	\end{pmatrix}
	\]
	is invertible, as its determinant takes the constant value $- \frac12$ irregardless of the profile $\zeta$; this is also reflected in the fact that $L$ and $\ub{L}$ are never collinear.
\end{rmk}

\section{Perturbed system}\label{sect:persys}
Our goal is to demonstrate that the simple wave solutions  constructed in the previous section are nonlinearly stable under the flow of \eqref{eq:mainwave} for sufficiently small perturbations. Note that for the trivial case $\zeta \equiv 0$ we recover small data global existence as a result.  

\begin{assm}\label{assm:zeta:hyp}
To ensure that corresponding initial value problem (with perturbed data prescribed at $t = 0$) is locally well-posed, we will make the assumption that 
\begin{equation}
	\inf_{x\in\Real} H(0) (\zeta'(x))^2 > -1
\end{equation}
for the background simple wave solution; see \thref{rmk:hyp}.
\end{assm}

We should mention at this juncture how we intend to \emph{compare} the background solution and the perturbed solution. 
Let us adorn as before the background quantities with a ring, so the background solution is $\mathring{\phi}$ and the perturbed solution is $\phi$. 
As both $\mathring{\phi}$ and $\phi$ are solutions to \eqref{eq:mainwave}, and hence are both functions defined on $(t,x)\in \Real^{1,1}$, the natural inclination is to compare the pointwise values $\mathring{\phi}(t,x) - \phi(t,x)$, and so forth for the derivatives. 
For our argument it is however more convenient to compare the two solutions via the conformal structure defined by the acoustic metric, or, more precisely, through the dynamical double null coordinates $u, \ub{u}$.

In particular, we will think of the background quantities $\mathring{\Psi}=0$, $\ub{\mathring{\Psi}}$, $\mathring{\sigma}=0$, $\mathring{\phi}$ as functions, not of the physical spacetime domain, but of some double null coordinate system $u, \ub{u}$. 
Similarly, we identify the perturbed quantities $\Psi, \ub{\Psi}, \sigma, \phi$ as functions of the \emph{same} dynamical coordinate system $u, \ub{u}$. 
Corresponding to these two solutions we reconstruct two transition mappings $\mathring{U}, U$ representing how, in each setting, the dynamical coordinates $u, \ub{u}$ are to be regarded as functions of the spacetime coordinates $(t,x)$. 
We will fix the gauge by requiring that along the initial data surface $\{t = 0\}$ the mappings $U$ and $\mathring{U}$ are identified. 
Therefore effectively we will be comparing, e.g. $\mathring{\phi}(t,x)$ with $\phi \circ U^{-1} \circ \mathring{U}(t,x)$. 

Taking this point of view has a couple advantages. 
Firstly, as seen already earlier in this paper, the equations of motion take particularly simple form in the dynamical coordinates $u, \ub{u}$, and simplifies the analysis. 
Secondly, this formulation makes it easier to factor out the effect of \emph{modified scattering} from the analysis; this is particularly convenient as in one spatial dimension solutions to the wave equation do not, in general, decay.

\begin{rmk}\label{rmk:phaseshift}
	To illustrate modified scattering, consider two simple wave solutions corresponding to $\phi(t,x) = \zeta(t-x)$ and $\tilde{\phi}(t,x) = \tilde{\zeta}(t-x)$, where $\zeta - \tilde{\zeta}$ is a small, compactly supported function. 
	For all points $t, x$ such that $\abs{t-x}$ is sufficiently large, we see then that the corresponding vector fields $L, \ub{L}$ and $\tilde{L}$ and $\ub{\tilde{L}}$ are equal (and all are in fact locally constant). 
	However, if we integrate the vector field $L$ and $\tilde{L}$ to obtain the level sets of the function $u$ and $\tilde{u}$, by requiring that they agree as $t-x \to -\infty$, we see that as $t-x \to +\infty$ there would (generally) be a phase-shift. 
\end{rmk}

Returning to the problem at hand, we will study all our \emph{wave} quantities as functions of $u, \ub{u}$. 
To emphasize the perturbative aspect of our analysis, we will write (regarding all functions as functions of $u, \ub{u}$, and not as the ambient coordinate $t,x$)
\begin{align}
	\psi & = \Psi - \mathring{\Psi} = \Psi;\\
	\ub{\psi} &= \ub{\Psi} - \ub{\mathring{\Psi}};
\end{align}
We do not introduce new notation for $\sigma$ as $\mathring{\sigma} = 0$.
The function $\zeta = \zeta(u)$, as before, denotes the profile of the simple wave background; we will write the perturbed solution as
\begin{equation}
	\xi = \phi - \mathring{\phi} = \phi - \zeta.
\end{equation}
A direct computation using \eqref{eq:semisys} gives us the following perturbed system for the ``wave variables'':
\begin{subequations}\label{eq:semisys:p}
	\begin{gather}
		\paru\parub \psi + \frac12 G(\sigma) \paru \sigma \parub \psi + \frac12 G(\sigma) \paru \psi \parub \sigma = 0, \label{eq:semisys:p1} \\
		\paru\parub \ub{\psi} + G(\sigma) \paru\sigma \zeta''(\ub{u}) + \frac12 G(\sigma) \paru \sigma \parub \ub{\psi} + \frac12 G(\sigma) \paru \ub{\psi} \parub \sigma = 0, \label{eq:semisys:p2} \\
		\sigma = - \psi (2 \zeta'(\ub{u}) + \ub{\psi}), \label{eq:semisys:p3}\\
		\paru\parub \xi + \frac14 \sigma( 1 + 2 f'(\sigma) \sigma) H'(\sigma) \left[ \paru \sigma \parub \xi + \paru \xi \parub\sigma + \zeta'(\ub{u}) \paru\sigma \right] = 0. \label{eq:semisys:p4}
	\end{gather}
\end{subequations}

\section{Statement of the main theorem}
For convenience we introduce the notation
\begin{equation}
	\Osub\epsilon\gamma{x} = \frac{\epsilon}{(1 + |x|)^{1+\gamma}}.
\end{equation}
Recall that a function $\phi$ of one variable $x$ is said to be \emph{of moderate decrease} if there exists some $M,\gamma > 0$ such that 
\[ \abs{\phi(x)} \leq \Osub{M}\gamma{x}. \]
In particular functions of moderate decrease are absolutely integrable. 

\begin{thm}\label{thm:main}
	Consider the initial value problem for \eqref{eq:mainwave} with initial data given in the rectangular coordinate system
	\[ \phi(0,x) = \phi_0(x), \qquad \partial_t \phi(0,x) = \phi_1(x), \]
	where $(\phi_0, \phi_1) \in C^2 \times C^1$. 
	Let $\zeta$ be a $C^2$ profile satisfying \thref{assm:zeta:hyp}, with some $M_\zeta, \bar{\gamma} > 0$ such that 
	\[ \abs{\zeta(x)}, \abs{\zeta'(x)}, \abs{\zeta''(x)} \leq \Osub{M_\zeta}{\bar{\gamma}}{x}. \]
	There exists a constant $\bar{\epsilon}$ depending on the value of $M_0$ of \thref{assm:f:assm}, and the values $M_\zeta$ and $\bar{\gamma}$, such that whenever the initial data satisfies 
\[ \abs{ \phi_0(x) - \zeta(-x) }, \abs{\phi_0'(x) + \zeta'(-x)}, \abs{\phi_0''(x) - \zeta''(-x)} , \abs{\phi_1(x) - \zeta'(-x)}, \abs{\phi_1'(x) + \zeta''(-x)} \leq \Osub{\bar{\epsilon}}{\bar{\gamma}}{x},\]

a unique global $C^2$ solution exists, with the solution depending Lipshitz-continuously on the initial data, when measured \emph{with respect to their respective dynamical double null coordinate systems}.  
\end{thm}

\section{Initial data and gauge fixing}

In this section we discuss the construction of the initial data for the reduced, semilinear system \eqref{eq:semisys:p}, and its smallness properties relative to prescribed smallness of the initial data $\phi$ in the rectangular coordinate system $(t,x)$. We will follow the notations introduced above and use ringed variables $\mathring{\phi}$ etc.\ to denote the variables associated to the background simple wave solution $\mathring{\phi}(t,x) = \zeta(t-x)$. 

The prescribed smallness condition immediately implies that, relative to the rectangular coordinate system, there is a universal constant $C$ such that 
\[\abs{\Psi(0,x) - \mathring{\Psi}(0,x)}, \abs{\ub{\Psi}(0,x) - \ub{\mathring{\Psi}}(0,x)}, \abs{\partial_x \Psi(0,x) - \partial_x \mathring{\Psi}(0,x)}, \abs{\partial_x \ub{\Psi}(0,x) - \partial_x \ub{\mathring{\Psi}}(0,x)} \leq \Osub{C\bar{\epsilon}}{\bar{\gamma}}{x}. \]
Similarly, there is a constant $C$ depending on $M_\zeta$ such that
\[  \abs{\sigma(0,x)}, \abs{\partial_x \sigma(0,x)} \leq \Osub{C\bar{\epsilon}}{1 + 2 \bar{\gamma}}{x}.\]
By the definitions \eqref{eq:dfn:acinvmetric} and \eqref{eq:dfn:acmetric} of the acoustic metric and its inverse, we see that this means that there is a constant $C$ depending on $M_\zeta$ and $M_0$ such that 
\begin{equation}\label{eq:acmetric:idbnd}
\abs{g_{\mu\nu}(0,x) - \mathring{g}_{\mu\nu}(0,x)} , \abs{g^{\mu\nu}(0,x) - \mathring{g}^{\mu\nu}(0,x)} \leq \Osub{C\bar{\epsilon}}{1 + 2\bar{\gamma}}{x} .
\end{equation}
In view of \thref{assm:zeta:hyp}, we see that provided $\bar{\epsilon}$ is sufficiently small, we must have $g^{00}(0,x) < 0$ (and hence $\{t = 0\}$ is spacelike). 
Therefore, there exists some $\bar{\epsilon}_0$ such that \emph{provided} $\bar{\epsilon} < \bar{\epsilon}_0$, we can solve \eqref{eq:mainwaveH} for $\partial^2_{tt}\phi(0,x)$, and obtain the following estimate for some constant $C$ (depending on $\bar{\epsilon}_0$, $M_0$, and $M_\zeta$)
\begin{equation}\label{eq:ttbnd}
	\abs{\partial_t \Psi(0,x) - \partial_t \mathring{\Psi}(0,x)}, \abs{\partial_t \ub{\Psi}(0,x) - \partial_t \ub{\mathring{\Psi}}(0,x)} \leq \Osub{C\bar{\epsilon}}{\bar{\gamma}}{x}. 
\end{equation}
\begin{rmk}\label{rmk:zeta:ind:bnd}
	We claim, furthermore, when $f'(0) \leq 0$, the constant $C$ in \eqref{eq:ttbnd} can be taken to be \emph{independent of $M_\zeta$} (note that when $f'(0) > 0$ our \thref{assm:zeta:hyp} places a strict size-limit on $M_\zeta$, and in that case $M_\zeta$ cannot be taken to be arbitrarily large). 

	First let us examine the case of \eqref{eq:ttbnd}. This requires solving \eqref{eq:mainwaveH} at $t = 0$ to find the values of $\partial^2_{tt}\phi$. Here we make use of the fact that the background solution $\mathring{\phi}$ and the corresponding $\partial^2_{tt}\mathring{\phi}$ value verifies $\partial_t \mathring{\phi} = - \partial_x \mathring{\phi}$ and $\partial^2_{tt} \mathring{\phi} = \partial^2_{xx}\mathring{\phi} = - \partial^2_{tx} \mathring{\phi}$. In particular, this means the product 
	\[   \eta^{\mu\alpha} \partial_\alpha \mathring{\phi} \partial^2_{\mu\nu} \mathring{\phi} \equiv 0. \]
	With this we see that the perturbation satisfies a reduced equation
\[\big( \eta^{\mu\nu} + 2 f'(0) \eta^{\mu\alpha}\partial_\alpha \mathring{\phi} \eta^{\nu\beta} \partial_\beta \mathring{\phi} + O(\partial(\phi - \mathring{\phi}))\big) \partial^2_{\mu\nu} (\phi - \mathring{\phi}) 	= - 2 f'(\sigma) \eta^{\mu\alpha}\partial_\alpha(\phi - \mathring{\phi}) \eta^{\nu\beta} \partial_\beta(\phi - \mathring{\phi}) \partial^2_{\mu\nu} \mathring{\phi}. \]
	The term $O(\partial(\phi - \mathring{\phi}))$ may have a coefficient depending on $M_\zeta$, but by choosing $\bar{\epsilon}$ sufficiently small this dependence can be overwhelmed, and is never an issue. 
	Attention, however, should be paid to what appears on the right of the equality sign: for generic perturbations of a linear system, we expect the right hand side depend only \emph{linearly} on $\phi - \mathring{\phi}$. The structure of the simple wave solutions means that here we have a quadratic dependence, and hence the potential largeness of the background $\partial^2_{\mu\nu}\mathring{\phi}$ can again be overwhelmed by sufficiently small $\bar{\epsilon}$. 

	Our final concern is with the linear part of the equation whose coefficients are $\eta^{\mu\nu} + 2 f'(0) \partial_{\alpha}\mathring{\phi} \partial_\beta \mathring{\phi} \eta^{\mu\alpha} \eta^{\nu\beta}$. Generically this depends on the background and may force $\partial^2_{tt}(\phi - \mathring{\phi})$ to be much larger than $\bar{\epsilon}$. 
	If we expand this in terms of $\zeta$, we find that the coefficients, expressed in matrix form, is
	\[ \begin{pmatrix} 
			-1 + 2 f'(0) (\zeta')^2 &  2 f'(0) (\zeta')^2 \\
			2 f'(0) (\zeta')^2 & 1 + 2f'(0) (\zeta')^2 
	\end{pmatrix} \]
	which means that when $f'(0) \leq 0$, the linear part of the solution bounds $\partial^2_{tt} (\phi(0,x) - \mathring{\phi}(0,x))$ in terms of $\partial^2_{tx}(\phi(0,x) - \mathring{\phi}(0,x))$ and $\partial^2_{xx} (\phi(0,x) - \mathring{\phi}(0,x))$ \emph{with a universal bound independent of either $M_0$ or $M_\zeta$}. 	
\end{rmk}

We next use our gauge freedom to set the values of $u$ and $\ub{u}$ at $t = 0$, by letting them satisfy
\begin{equation}
	u(0,x) = x, \qquad \ub{u}(0,x) = -x,
\end{equation}
and requiring $\partial_t \ub{u}(,x) > 0$ (this latter can be achieved since for the planewave background one has $\partial_t\ub{\mathring{u}}= 1$; the value of $\partial_t u$ will then be uniquely specified by requiring the two are independent). 
By the eikonal equations, the value of the time derivatives $\partial_t u(0,x)$ and $\partial_t \ub{u}(0,x)$ initially can be found by solving pointwise the quadratic equations
\[  (\partial_t u)^2  - 2 f'(\sigma) ( - \partial_t u \cdot \partial_t\phi + \partial_x \phi)^2 = 1 =  (\partial_t \ub{u})^2 - 2 f'(\sigma) (  \partial_t \ub{u} \cdot \partial_t\phi + \partial_x \phi)^2; \]
notice that the solution is $C^1$ in $x$. 
We further have, for some $C$ depending on $M_0$, and $M_\zeta$, 
\begin{equation}\label{eq:ut:idbnd}
\abs{\partial_t u(0,x) - \partial_t \mathring{u}(0,x)}, \abs{\partial_t \ub{u}(0,x) - \partial_t \mathring{\ub{u}}(0,x)} \leq \Osub{C\bar{\epsilon}}{\bar{\gamma}}{x}.
\end{equation}

Now, considering the variables $\psi, \ub\psi, \sigma, \xi$ \emph{as functions of the dynamical coordinates $u, \ub{u}$} (as described in Section \ref{sect:persys}), and noting that the curve $\{t = 0\}$ is now the curve $\{u = - \ub{u}\}$ by construction, we find:
\begin{lem}
	For any $\epsilon_0 > 0$, we can choose $\bar{\epsilon}$ sufficiently small under the hypotheses of \thref{thm:main}, such that the initial conditions along $\{ u = - \ub{u} \}$ satisfy
	\begin{multline}\label{eq:semisys:datasmall}
		\abs{\psi(s,-s)}, \abs{\paru\psi(s,-s)}, \abs{\parub\psi(s,-s)}, 
		\abs{\ub{\psi}(s,-s)}, \abs{\paru\ub{\psi}(s,-s)}, \\
		\abs{\parub\ub{\psi}(s,-s)}, 
		\abs{\sigma(s,-s)}, \abs{\xi(s,-s)}, \abs{\paru\xi(s,-s)}, \abs{\parub\xi(s,-s)} \leq \Osub{\epsilon_0}{\bar{\gamma}}{s}. 
	\end{multline}
\end{lem}
	
\section{Global existence of the semilinear system}
In this section we prove the following proposition. 

\begin{prop}\label{prop:semi:gwp}
	Consider the initial value problem to the semilinear system \eqref{eq:semisys:p} with initial data given on $\{u + \ub{u} = 0\}$. There exists a $\delta_0 > 0$ such that for every $\delta \in (0, \delta_0)$, there exists an $\epsilon' > 0$ such that if the initial data satisfy \eqref{eq:semisys:datasmall} with $\epsilon_0 < \epsilon'$, then the solution exists globally (in the $u, \ub{u}$ coordinates) with bounds
	\begin{gather*}
		\abs{\psi}, \abs{\ub{\psi}}, \abs{\xi} \leq \delta, \\
		\abs{\paru\psi}, \abs{\paru\ub{\psi}}, \abs{\paru\xi}  \leq \Osub{\delta}{\bar{\gamma}}{u}, \\
		\abs{\parub\psi}, \abs{\parub\ub{\psi}}, \abs{\parub\xi} \leq \Osub{\delta}{\bar{\gamma}}{\ub{u}},
	\end{gather*}
\end{prop}

\begin{rmk}
	Here we carry out our estimates by directly estimating the $L^\infty$ norms using the fundamental solution of the wave equation for the solution itself, and using the transport equations for the derivatives. One can alternatively approach the same problem using weighted $L^2$ energy estimates; see \cite{LuYaYu2017}. 
\end{rmk}

The following basic calculus result is convenient.
\begin{lem} \label{lem:calc}
	$\int_{\Real} \Osub{\epsilon}{\gamma}{s} \D{s} \leq 2\epsilon (1 + \gamma^{-1})$.
\end{lem}

We first consider the closed system \eqref{eq:semisys:p1}--\eqref{eq:semisys:p3}.
Given $(\tilde{\psi}, \ub{\tilde{\psi}}) \in C^1(\Real^2; \Real^2)$, we can set $\tilde{\sigma}\in C^1(\Real^2)$ in accordance to the algebraic \eqref{eq:semisys:p3}. (Note here we consider the domain $\Real^2$ as the set $\{(u, \ub{u})\}$.)
For the first two equations, observe that the second term in \eqref{eq:semisys:p2} contains terms that are \emph{linear} in the perturbation variables. 
This means that for the contraction mapping argument that we will use, some care needs to be made in regards to that term. 
Here we will use the fact that we have a \emph{nilpotent} structure: the linear term is essentially of the form 
\[  G(0) \paru\psi \zeta'(\ub{u}) \zeta''(\ub{u}) \]
and contributes to the equation for $\ub{\psi}$, so it is ``off diagonal''. Hence we can proceed by first solving for $\psi$ and then using it as the input for $\ub{\psi}$. 

\subsection{Iteration scheme}

First, let us denote in the sequel the function
\begin{equation}
	\sigma(\psi, \ub{\psi}) \eqdef - \psi (2 \zeta'(\ub{u}) + \ub{\psi}).
\end{equation}
Consider the operator $T$ which sends $(\tilde{\psi}, \ub{\tilde{\psi}}) \in C^1(\Real^2; \Real^2)$ to the solution $(\psi, \ub{\psi})$ where we first solve
\[ \paru\parub \psi = - \frac12 G(\sigma(\tilde{\psi}, \ub{\tilde{\psi}})) \paru  \sigma(\tilde{\psi}, \ub{\tilde{\psi}})\parub \tilde{\psi} - \frac12 G(\sigma(\tilde{\psi}, \ub{\tilde{\psi}})) \paru \tilde\psi \parub \sigma(\tilde{\psi}, \ub{\tilde{\psi}}).\]
Having solved for $\psi$, we next solve
\[\paru\parub \ub{\psi} = -G(\sigma(\psi, \ub{\tilde{\psi}})) \paru \sigma(\psi, \ub{\tilde{\psi}})\zeta''(u) - \frac12 G(\sigma(\psi, \ub{\tilde{\psi}})) \paru \sigma(\psi, \ub{\tilde{\psi}}) \parub \ub{\tilde\psi} - \frac12 G(\sigma(\psi, \ub{\tilde{\psi}})) \paru \ub{\tilde\psi} \parub \sigma(\psi, \ub{\tilde{\psi}}).\]
It is clear that $T:C^1(\Real^2; \Real^2) \to C^1(\Real^2; \Real^2)$ given our assumptions, by re-writing the equation as inhomogeneous wave equations in integral form.  

Denote by $\mathfrak{X}_\epsilon$ the subset
\begin{multline}
	\mathfrak{X}_\epsilon \eqdef \bigl\{ (\psi, \ub{\psi}) \in C^1(\Real^2; \Real^2) : 
		\abs{\psi}\leq \epsilon^2 , \abs{\ub{\psi}} \leq \epsilon, 
	\abs{\paru\psi} \leq \Osub{\epsilon^2}{\bar\gamma}{u}, \\
	\abs{\paru\ub{\psi}} \leq \Osub{\epsilon}{\bar{\gamma}}{u},
\abs{\parub\psi}\leq \Osub{\epsilon^2}{\bar{\gamma}}{\ub{u}}, \abs{\parub\ub{\psi}} \leq \Osub{\epsilon}{\bar{\gamma}}{\ub{u}} \bigr\}.
\end{multline}
We will show that when $\epsilon_0, \delta$ are sufficiently small, the mapping $T: \mathfrak{X}_\delta \to \mathfrak{X}_\delta$ is a contraction mapping; this argument is largely standard, except for the use of the nilpotent structure above. 

\subsection{Range of $T$}
First we check that for $\epsilon_0, \delta$ sufficiently small, $T$ maps $\mathfrak{X}_\delta$ into itself. 
First observe that for $\delta$ sufficiently small compared to $M_\zeta$, we have that 
\[ \sigma( \mathfrak{X}_\delta ) \subset [-m_0, m_0] \]
(where $m_0$ is of the \thref{assm:f:assm}). And in particular $G\circ \sigma$ and $G'\circ\sigma$ are bounded by $M_0$ on $\mathfrak{X}_\delta$.  

The derivatives of $\sigma$ can be bounded by the product rule, the obvious bounds on $\mathfrak{X}_\delta$ are
\[ \abs{\paru\sigma} \leq \Osub{3M_\zeta \delta^2}{\bar{\gamma}}{u}, \quad \abs{\parub\tilde{\sigma}} \leq \Osub{4M_\zeta \delta^2}{\bar{\gamma}}{\ub{u}}. \]
(In the latter we need the moderate decrease of the background simple wave.) 
The linear wave equation satisfied by $\psi$ then has the form
\[ \abs{\paru\parub \psi} \leq 4 M_0 M_\zeta \Osub{\delta^2}{\bar{\gamma}}{u} \Osub{\delta^2}{\bar{\gamma}}{\ub{u}} .\]
We can integrate in $u$ to find, for example
\[ \parub\psi(u, \ub{u}) = \int_{-\ub{u}}^{u} \paru\parub\psi(s, \ub{u}) \D{s} + \parub\psi(-\ub{u}, \ub{u}) \]
which yields, in view of the initial data bounds and Lemma \ref{lem:calc},
\[ \abs{\parub\psi(u, \ub{u})} \leq \Osub{\epsilon_0}{\bar{\gamma}}{\ub{u}} + 8 M_0 M_\zeta (1 + (\bar\gamma)^{-1}) \delta^2 \Osub{\delta^2}{\bar{\gamma}}{\ub{u}}. \]
Integrating this in $\ub{u}$ we further obtain 
\[ \abs{\psi} \leq \epsilon_0 + 2 \epsilon_0 (1 + (\bar\gamma)^{-1}) + 16 M_0 M_\zeta \delta^4 (1 + (\bar\gamma)^{-1})^2.\]
This shows clearly that for $\delta$ sufficiently small, the pair $(\psi, \ub{\tilde{\psi}})$ \emph{is also an element of $\mathfrak{X}_{\delta}$}.

The linear wave equation satisfied by $\ub{\psi}$ on the other hand looks like (now using $(\psi, \ub{\tilde{\psi}})$ as the input to $\sigma$)
\[ \abs{\paru\parub \ub{\psi}} \leq 4 M_0 M_\zeta( 1 + \delta) \Osub{\delta}{\bar{\gamma}}{u} \Osub{\delta}{\bar{\gamma}}{\ub{u}}\]
with the term with coefficient $1$ coming from the linear dependence on $\sigma$ which depends linearly on $\psi$ and so the corresponding term has only $\delta^2$ smallness.  
Integrating we get
\[ \abs{\parub\ub{\psi}(u, \ub{u})} \leq \Osub{\epsilon_0}{\bar{\gamma}}{\ub{u}} + 12 M_0 M_\zeta (1 + (\bar{\gamma})^{-1}) \delta \Osub{\delta}{\bar{\gamma}}{\ub{u}} \]
and
\[ \abs{\ub\psi} \leq \epsilon_0 + 2 \epsilon_0 (1 + (\bar\gamma)^{-1}) + 24 M_0 M_\zeta \delta^2 (1 + (\bar{\gamma})^{-1})^{2}. \] 
So provided
\begin{equation}\label{eq:epsilonsmallness} 6(1 + (\bar{\gamma})^{-1}) \epsilon_0 \leq \delta^2, \qquad \delta \leq \frac{1}{48 M_0 M_\zeta (1 + (\bar{\gamma})^{-1})^2} 
\end{equation}
we see that the whole argument goes through and the solution $(\psi, \ub{\tilde{\psi}}) \in \mathfrak{X}_\delta$. 

\subsection{Contraction}
We will prove $T$ is a contraction with respect to the metric
\begin{multline}
	d((\psi_1, \ub{\psi}_1), (\psi_2, \ub{\psi}_2)) = \max\Bigl\{ \norm[L^\infty]{\psi_1 - \psi_2}, \norm[L^\infty]{\ub{\psi}_1 - \ub{\psi}_2}, \\
		\norm[L^\infty]{(1 + \abs{u})^{1+\bar\gamma}(\paru\psi_1 - \paru\psi_2)}, \norm[L^\infty]{(1 + \abs{u})^{1+\bar\gamma}(\paru\ub{\psi}_1 - \paru\ub{\psi}_2)},\\
		\norm[L^\infty]{(1 + \abs{\ub{u}})^{1+\bar\gamma}(\parub\psi_1 - \parub\psi_2)}, 
	\norm[L^\infty]{(1 + \abs{\ub{u}})^{1+\bar\gamma}(\parub\ub{\psi}_1 - \parub\ub{\psi}_2)} \Bigr\},
\end{multline}
provided $\delta$ is sufficiently small (with $\epsilon_0$ adjusted suitably in accordance to \eqref{eq:epsilonsmallness})

Notice that on $\mathfrak{X}_\delta$ we have that the polynomial function $\sigma$ has bounded derivatives
\[  \abs{\frac{\partial}{\partial \psi} \sigma} \leq 3 M_\zeta, \quad \abs{\frac{\partial}{\partial\ub{\psi}}\sigma} \leq \delta^2.\]

Now, given $(\tilde{\psi}_1, \ub{\tilde{\psi}}_1)$ and $(\tilde{\psi}_2, \ub{\tilde{\psi}}_2)$, we can solve first for the corresponding $\psi_1, \psi_2$. Their difference solves a wave equation with vanishing initial data, whose inhomogeneity can be thus bounded by (up to a numerical constant $C$ which we won't make precise)
\[  C M_0 M_\zeta \Osub{\delta}{\bar{\gamma}}{u} \Osub{\delta}{\bar{\gamma}}{\ub{u}} d((\tilde{\psi}_1, \ub{\tilde{\psi}}_1), (\tilde{\psi}_2, \ub{\tilde{\psi}}_2)). \]
Therefore integrating as above we have that 
\begin{equation}\label{eq:psicontract}
d((\psi_1, 0), (\psi_2,0)) \leq  C M_0 M_\zeta \delta^2 d((\tilde{\psi}_1, \ub{\tilde{\psi}}_1), (\tilde{\psi}_2, \ub{\tilde{\psi}}_2)). 
\end{equation}

One can treat the \emph{nonlinear} contributions to the difference of $\ub{\psi}_1 - \ub{\psi}_2$ similarly, with the factor of $\delta^2$ above replaced by $\delta$, and we will not dwell on them. 
For the linear term, however, we must make use of the nilpotent structure.
To wit, the best estimate available for the difference
\[ G(\sigma(\tilde{\psi}_1, \ub{\tilde{\psi}}_1)) \paru \sigma(\tilde{\psi}_1, \ub{\tilde{\psi}}_1)\zeta''(\ub u) - G(\sigma(\tilde{\psi}_2, \ub{\tilde{\psi}}_2)) \paru \sigma(\tilde{\psi}_2, \ub{\tilde{\psi}}_2)\zeta''(\ub u) \]
would be of size 
\[  M_\zeta^2 M_0  d((\tilde{\psi}_1, \ub{\tilde{\psi}}_1), (\tilde{\psi}_2, \ub{\tilde{\psi}}_2)) \]
which will fail to give a contraction mapping for large $M_\zeta$. 
However, our use of the nilpotent structure allows us to estimate with \emph{the already solved} $\psi_1, \psi_2$ instead
\[ G(\sigma({\psi}_1, \ub{\tilde{\psi}}_1)) \paru \sigma({\psi}_1, \ub{\tilde{\psi}}_1)\zeta''(\ub u) - G(\sigma({\psi}_2, \ub{\tilde{\psi}}_2)) \paru \sigma({\psi}_2, \ub{\tilde{\psi}}_2)\zeta''(\ub u) \]
which can be bounded by 
\[ M_\zeta^2 M_0 d((\psi_1,0), (\psi_2,0)) + M_0 M_\zeta^2 \delta^2 d((0, \tilde{\psi}_1), (0, \tilde{\psi}_2)). \]
Bounding the first term by \eqref{eq:psicontract} which we have already obtained, we can in fact conclude that there exists some numerical constant $C$ for which 
\begin{equation}
	d((0,\ub{\psi}_1), (0,\ub{\psi}_2)) \leq  C M_0^2 M_\zeta^3 \delta^2 d((\tilde{\psi}_1, \ub{\tilde{\psi}}_1), (\tilde{\psi}_2, \ub{\tilde{\psi}}_2)). 
\end{equation}
And hence by further shrinking $\delta$ we can guarantee that $T$ is a contraction mapping.  

\subsection{Equation for $\xi$}
Having solved for $\sigma$, the existence and uniqueness, and estimates for $\xi$ follows simply from the \emph{linear} equation \eqref{eq:semisys:p4}. Here we only need to observe that that equation can be schematically written as
\[ \paru\parub \xi = A \parub \xi + \ub{A} \paru\xi + B \]
where the coefficients have the bounds
\begin{gather*}
	A \leq \frac{10 M_\zeta^2 M_0^2 \delta^2}{(1 + \abs{u})^{1 + \bar{\gamma}}},\\
	\ub{A} \leq \frac{10 M_\zeta^2 M_0^2 \delta^2}{(1 + \abs{\ub{u}})^{1 + \bar{\gamma}}},\\
	B \leq \frac{10 M_\zeta^3 M_0^2 \delta^2}{(1 + \abs{\ub{u}})^{1 + \bar\gamma}(1 + \abs{u})^{1 + \bar\gamma}}.
\end{gather*}
Directly integrating, using that the coefficients are quadratic in $\delta$, we see that if $\delta$ is sufficiently small we can find a unique solution satisfying the requisite bounds. 

\section{Jacobian bounds} 
In this section we study the geodesic equations \eqref{eq:geodesic}. 
Observing that
\[ 
	g^{\mu\nu} \Phi_\nu = \ub{L}^\mu \parub\phi + L^\mu \paru\phi 
\]
and 
\[
	L^\beta \Phi_\beta = \Omega^{-1} \parub \phi, \qquad \ub{L}^\beta \Phi_\beta = \Omega^{-1} \paru\phi,
\]
we rewrite in the transport equations in the following form:
\begin{equation}\label{eq:geod2}
	\begin{gathered}
		\begin{split}
			\parub L^\mu  &+ L^\mu \paru\phi H(\sigma) L^\gamma \parub \Phi_\gamma + \ub{L}^\mu \parub \phi H(\sigma) L^\gamma \parub \Phi_\gamma \\
				     &+ \Omega^{-1} H'(\sigma) \parub\phi \left[ L^\mu \paru \phi \parub\sigma + \frac12 \ub{L}^\mu \parub\sigma\parub\phi - \frac12 L^\mu \paru\sigma \parub\phi \right] = 0,
		\end{split} \\
		\begin{split}
			\paru \ub{L}^\mu &+ \ub{L}^\mu \parub\phi H(\sigma) \ub{L}^\gamma \paru \Phi_\gamma + L^\mu \paru \phi H(\sigma) \ub{L}^\gamma \paru \Phi_\gamma \\
					 &+ \Omega^{-1} H'(\sigma) \paru\phi \left[ \ub{L}^\mu \parub \phi \paru\sigma + \frac12 L^\mu \paru\sigma\paru\phi - \frac12 \ub{L}^\mu \parub\sigma \paru\phi \right] = 0.
		\end{split}
	\end{gathered}
\end{equation}
The key to our argument is the fact that the inhomogeneous terms in the transport equation do \emph{not} contain dangerous non-decaying terms: the equation for $\parub L^\mu$ sees every term in the inhomogeneity containing a $\parub$ derivative of $\Psi, \ub{\Psi}, \phi,$ or $\sigma$; while the equation for $\paru \ub{L}^\mu$ sees every term in the inhomogeneity containing a $\paru$ derivative. 
\emph{This is a consequence of the null condition enjoyed by our system.} 
This is in contrast with the setting with genuine nonlinearity, where one expects in the inhomogeneity a term with non-decaying coefficients and purely transversal derivatives. 

\begin{rmk}
	One can compare \eqref{eq:geod2} to equation (1.3.1) in \cite{SpHoLW2016}. There the quantity $\upmu$ is roughly $(L^0)^{-1}$. The problematic inhomogeneity $\frac12 G_{LL} \breve{X}\Psi$ should be compared to the first term in our inhomogeneity of $\parub L^\mu$, namely the term $L^\mu \paru\phi H(\sigma) L^\gamma \parub\Phi_{\gamma}$. The difference is that in the genuinely nonlinear case we lack the compensating factor of $\parub\Phi_{\gamma}$, in which case this term will give a Riccati type term that allows driving $L^\mu$ to $0$ in finite time, causing shock formation. 
\end{rmk}

\begin{prop}\label{prop:geom:gwp}
	There exists a $\varepsilon_0 > 0$ such that for every $\varepsilon \in (0,\varepsilon_0)$, there are corresponding values of $\delta', \epsilon''$, such that if $\psi, \ub{\psi}, \xi$ satisfy the bounds of \thref{prop:semi:gwp} with $\delta < \delta'$, and that the initial data in \thref{thm:main} satisfies $\bar{\epsilon} < \epsilon''$, the system \eqref{eq:geod2} can be solved for all $(u, \ub{u})$ with the uniform bound 
	\[  \abs{L^\mu - \mathring{L}^\mu}, \abs{\ub{L}^\mu - \ub{\mathring{L}}^\mu} < \varepsilon;\]
	here $\mathring{L}^\mu$ and $\ub{\mathring{L}}^\mu$ are corresponding values of the simple wave background, given by \eqref{eq:Lbg}. 
\end{prop}

Let us first derive some preliminary estimates using the bounds derived for \thref{prop:semi:gwp}. 
First, we have that by our setup the scalars
\begin{equation}
	\abs{\Psi} = \abs{\psi} \leq \delta, \qquad \abs{\ub{\Psi}} = \abs{\ub{\psi} + 2 \zeta'(\ub{u})} \leq \delta + \Osub{M_\zeta}{\bar{\gamma}}{\ub{u}}.
\end{equation}
This decomposition also tells us that
\begin{equation}
	\abs{\parub \Phi_\mu} \leq \Osub{2 M_\zeta}{\bar{\gamma}}{\ub{u}}, \quad \abs{\paru \Phi_\mu} \leq \Osub{2 \delta}{\bar{\gamma}}{u}. 
\end{equation}
Similarly, writing $\phi = \zeta(\ub{u}) + \xi$ we have 
\begin{equation}
	\abs{\parub \phi} \leq \Osub{2 M_\zeta}{\bar{\gamma}}{\ub{u}}, \quad \abs{\paru\phi} \leq \Osub{\delta}{\bar{\gamma}}{u}.
\end{equation}
The acoustic metric can be written, by way of \eqref{eq:dfn:acmetric}, as
\[
	g = - \D*{t}^2 + \D*{x}^2 - \frac{f'(\sigma)}{2 + 4 f'(\sigma)\sigma} \left[ \Psi^2 (\D*{t} + \D*{x})^2 + \ub{\Psi}^2 (\D*{t} - \D*{x})^2 + 2 \sigma (- \D*{t}^2 + \D*{x}^2) \right]. 
\]
Denote by $\mathring{g}$ the acoustic metric for the background simple wave solution, we have that there exists a constant $C$ depending only on $M_0$ and $M_\zeta$, such that the rectangular components 
\begin{equation}\label{eq:acmetric:fbnd}
	\abs{g_{\mu\nu} - \mathring{g}_{\mu\nu}} , \abs{g^{\mu\nu} - \mathring{g}^{\mu\nu}} \leq C \delta. 
\end{equation}
And hence the rectangular components of $g$ is arbitrarily close to that of $\mathring{g}$, provided we choose $\delta$ small enough in \thref{prop:semi:gwp}. 
In particular, we have that 
\[ \abs{\Omega^{-1} } \leq C ( \abs{L} \cdot \abs{\ub L} ). \]

From this we see that \emph{most} of the terms of the inhomogeneity in \eqref{eq:geod2} are harmless: in fact we can rewrite \eqref{eq:geod2} schematically as (for some constant $C$ depending only on $M_\zeta, M_0$)
\begin{align*}
	\parub L^\mu + \ub{L}^\mu \parub \phi H(\sigma) L^\gamma \parub \Phi_{\gamma} &= \frac{C\delta }{(1 + \abs{\ub{u}})^{1 + \bar{\gamma}}} (\abs{L}^2 + \abs{L}^4 + \abs{\ub{L}}^2 + \abs{\ub{L}}^4), \\
	\paru \ub{L}^\mu &= \frac{C\delta }{(1 + \abs{u})^{1 + \bar{\gamma}}} (\abs{L}^2 + \abs{L}^4 + \abs{\ub{L}}^2 +      \abs{\ub{L}}^4).
\end{align*}
For the first equation, we can further expand $\phi = \xi + \zeta$ and similarly $\Phi_\gamma$ in terms of $\psi, \ub{\psi}$ and $\zeta'$, to obtain
\[ \parub L^\mu + \ub{L}^\mu \zeta'(\ub{u}) H(0) \left( L^0 - L^1 \right) \zeta''(\ub{u}) 
	= \frac{C\delta }{(1 + \abs{\ub{u}})^{1 + \bar{\gamma}}} (\abs{L}^2 + \abs{L}^4 + \abs{\ub{L}}^2 +      \abs{\ub{L}}^4).
\]
Due to the presence of the $\delta$ smallness, these terms appearing on the right of the equality sign are essentially harmless for the contraction mapping principle argument, provided we take $\delta$ sufficiently small. 
As the arguments relating to those terms are standard, for brevity we will omit them from consideration and examine instead the model system
\begin{equation}\label{eq:geod3mod}
	\begin{gathered}
		\paru \ub{L}^\mu = 0, \\
		\parub L^\mu + \zeta'\zeta'' H(0) \ub{L}^\mu (L^0 - L^1) = 0. 
	\end{gathered}
\end{equation}

We will treat \thref{prop:geom:gwp} perturbatively, comparing against the background solutions $\mathring{L}$ and $\ub{\mathring{L}}$. Writing
\begin{align*}
\ell^\mu \eqdef L^\mu - \mr L^\mu, \\
\ub\ell^\mu \eqdef \ub L^\mu - \ub{\mr L}^\mu,
\end{align*}
for the difference between the background and the perturbed solutions, we find that our model equation \eqref{eq:geod3mod} becomes
\begin{equation}\label{eq:geod4mod}
	\begin{gathered}
		\paru \ub{\ell}^\mu = 0, \\
		\parub \ell^\mu + \zeta'\zeta'' H(0) \left[(\ell^0 - \ell^1) \ub{\mathring{L}}^\mu +  \ub{\ell}^\mu (-2 + \ell^0 - \ell^1) \right]= 0. 
	\end{gathered}
\end{equation}
We note that in deriving \eqref{eq:geod4mod}, the explicit form of $\mr L, \ub{\mr L}$ led to a crucial cancellation of $2 H(0) \zeta' \zeta''$ which has no smallness factor. Moreover, we highlight the presence of terms \emph{linear} in $\ell$ in the second equation \emph{without} any small coefficients, so these terms cannot be controlled as nonlinear inhomogeneities: we must consider the corresponding linear evolution. We will derive smallness from the smallness of the initial data which we can prescribe. 

In particular, as a consequence of the bounds \eqref{eq:acmetric:idbnd} and \eqref{eq:ut:idbnd}, we have that there exists some constant $C$ (which may depend on $M_0, M_\zeta$), such that the components
\[ \abs{\ell^\mu(s,-s)}, \abs{\ub{\ell}^\mu(s,-s)} \leq \Osub{C\bar{\epsilon}}{\bar{\gamma}}{s}. \]

We observe now that the system \eqref{eq:geod4mod} is essentially linear: solving the first equation guarantees that $\ub{\ell}$ is globally bounded by the initial data bound $C\bar{\epsilon}$.
The second equation enjoys a further decomposition, using that $\ub{\mathring{L}} = (-1,-1)$: we have the decoupled linear equations for $\ell^0 \pm \ell^1$
\begin{equation} \label{eq:geo5mod}
\begin{gathered}
	\parub( \ell^0 - \ell^1) + \zeta' \zeta'' H(0) (\ub{\ell}^0 - \ub{\ell}^1) (-2 + \ell^0 - \ell^1) = 0,  \\
	\parub( \ell^0 + \ell^1) - 2 \zeta' \zeta'' H(0) (\ell^0 - \ell^1) + \zeta' \zeta'' H(0) (\ub{\ell}^0 + \ub{\ell}^1) (-2 + \ell^0 - \ell^1) = 0.
\end{gathered}
\end{equation}
Key here is that $\zeta'$ and $\zeta''$ are assumed to decay in $\ub{u}$; particular we have that $\zeta' \zeta'' \leq \Osub{M_\zeta^2}{1 + 2 \bar{\gamma}}{\ub{u}}$ and is integrable in $\ub{u}$. 
From this we see that there exists a constant $C$ depending on $M_0$ and $M_\zeta$ such that under the assumptions of \thref{thm:main}, the linearized system \eqref{eq:geo5mod} can be solved first for $\ell^0 - \ell^1$ with a uniform bound of $C\bar{\epsilon}$. Using this smallness as input for the equation for $\ell^0 + \ell^1$, in total this means that with $\ell, \ub{\ell}$ can be solved by a contraction mapping argument and are uniformly bounded by $C\bar{\epsilon}$. 
\section{Conclusion}

\thref{thm:main} now follows by combining \thref{prop:semi:gwp} and \thref{prop:geom:gwp}. 
In particular, the acoustic metric bound \eqref{eq:acmetric:fbnd} implies that for sufficiently small $\delta$ we can guarantee that the metric is non-degenerate, which combined with the result of \thref{prop:geom:gwp} shows that for sufficiently small $\bar{\epsilon}$ we can guarantee that $L^0$ and $\ub{L}^0$ are both bounded away from zero, and that $\Omega$ remains in a neighborhood of the background value $-\frac12$. And hence by the argument in \thref{rmk:cov:reg} we see that the transformation between the $(u, \ub{u})$ coordinate system and the $(t,x)$ system is a $C^1$ diffeomorphism of $\mathbb{R}^2$. 

Note finally that as $\Psi, \ub{\Psi}$ are $C^1$ functions of the coordinates $(u, \ub{u})$, and they represent the values of \emph{rectangular coordinate derivatives} of $\phi$, that $(u, \ub{u})$ is $C^1$ is sufficient to provide the reverse transformation to guarantee that the solution $\phi$ is $C^2$ measured with respect to the rectangular coordinate system. 

Uniqueness of the solution follows from the fact that the existence can be proven using a contraction mapping argument; the same also shows we have in fact \emph{Lipshitz} dependence of the solution on initial data. 
\begin{rmk}
	Generally for quasilinear wave equations one only expects continuous dependence of solution on the initial data; in our case Lipshitz dependence can be achieved because of the strong semilinearization of our equations when expressed in the double-null coordinates. 
\end{rmk}

\printbibliography
\end{document}